%% file: main.tex
\documentclass[runningheads]{llncs}
\usepackage{graphicx}
\usepackage{subfig}
\usepackage{caption}
\usepackage[utf8]{inputenc}
\usepackage{amsfonts, amsmath, amssymb}
\usepackage{xcolor}
\usepackage{booktabs}
\usepackage{tablefootnote}
\usepackage{array}
\usepackage{float} 
\usepackage{listings}
\newcolumntype{L}[1]{>{\raggedright\let\newline\\\arraybackslash\hspace{0pt}}m{#1}}
\newcolumntype{C}[1]{>{\centering\let\newline\\\arraybackslash\hspace{0pt}}m{#1}}
\newcolumntype{R}[1]{>{\raggedleft\let\newline\\\arraybackslash\hspace{0pt}}m{#1}}

\lstdefinelanguage{SCENIC}{
    keywords={param, VerifaiRange, behavior, try, interrupt, Uniform, at, with, require},
    keywordstyle={\color{blue}\bfseries},
    sensitive=false, 
    morecomment=[l]{\#}, 
    showstringspaces=false,
    numbers=left,
    stepnumber=1,
} %

\input macros

\begin{document}

\title{Specification-Guided Critical Scenario Identification for Automated Driving
\thanks{\includegraphics[height=0.25cm]{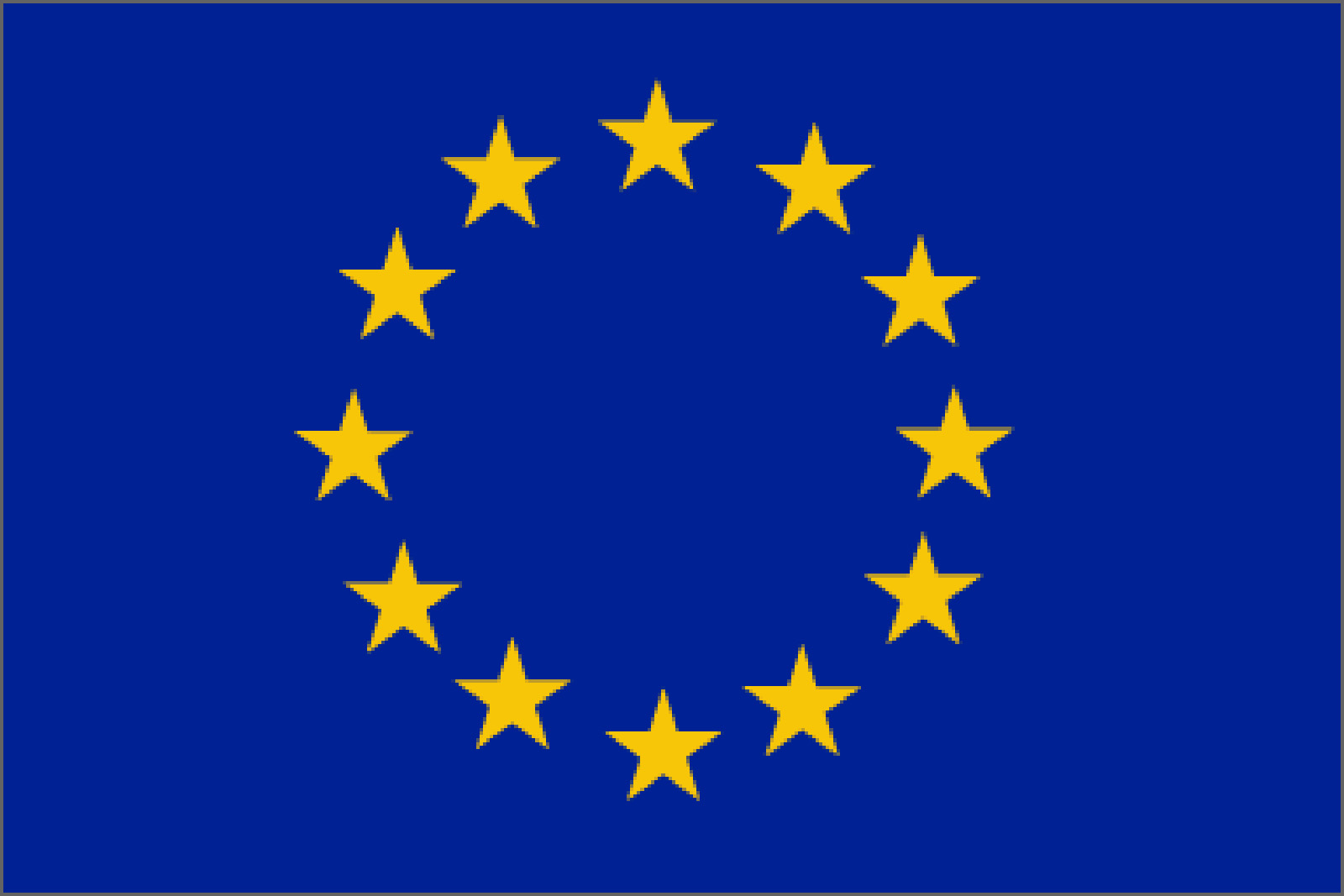} This project has received funding from the European Union’s Horizon 2020 research and innovation programme under grant agreement No 956123.}
}

\author{Adam Molin\inst{1} \and
Edgar A. Aguilar\inst{2} \and
Dejan Ni\v{c}kovi\'{c}\inst{2}\and
Mengjia Zhu\inst{3} \and
Alberto~Bemporad\inst{3} \and
Hasan Esen\inst{1}}

\authorrunning{A. Molin et al.}

\institute{DENSO AUTOMOTIVE Deutschland GmbH, 85386 Eching, Germany
\email{\{a.molin,h.esen\}@eu.denso.com} \and
AIT Austrian Institute of Technology GmbH, 1210 Vienna, Austria
\email{\{edgar.aguilar,dejan.nickovic\}@ait.ac.at}\\
 \and
IMT School for Advanced Studies Lucca, 55100 Lucca, Italy\\
\email{\{mengjia.zhu,alberto.bemporad\}@imtlucca.it}}

\maketitle

\begin{abstract}
To test automated driving systems, we present a case study for finding critical scenarios in driving environments guided by formal specifications. To that aim, we devise a framework for critical scenario identification, which we base on open-source libraries that combine scenario specification, testing, formal methods, and optimization. 

\end{abstract}

\section{Introduction}
With the complexity of the automated driving (AD) system and its driving environment, verification and validation (V\&V) is regarded as one of the major challenges of AD development \cite{zhang2022finding}. \emph{Scenario-based testing} (SBT) was introduced as an essential method for facilitating the overall safety assurance of ADs. In SBT, the expected behavior of an AD system is described by a representative set of scenarios that are relevant for its safe use. The SBT paradigm facilitates shifting the AD testing from the physical to the simulation environment. The use of virtual testing has manifold advantages -- more specifically it allows to: (1) explore efficiently a large number of situations originating from the catalog of relevant scenarios, (2) reproduce environment conditions (fog, night, rain, etc.) that are hard to enforce in a physical environment, and (3) play dangerous scenarios without risk to humans, other vehicles or infrastructure. 

Despite significant advances in research and standardization of SBT, there are still remaining open issues. One of them is to determine the critical scenarios among the virtually infinite number of scenarios with an abundance of influential factors ranging from weather or road conditions, to the behaviors of surrounding road users.
A first attempt to keep the number of scenarios manageable is to restrict the operational design domain (ODD) of the AD system. According to \cite{orad_j3016}, the ODD is defined as the operating conditions under which a given AD system is specifically designed to function. However, there are some factors, including the dynamic behavior of the road users, which cannot be controlled in the ODD. Thus, efficient methods are needed to identify the critical scenarios from the scenario space within the ODD. An extensive survey study on finding critical scenarios has been conducted in \cite{zhang2022finding}. With regard to specification-guided critical scenario identification, our work is closely related to \cite{DBLP:journals/jar/DreossiDS19,qin2019automatic,tuncali2019requirements,DBLP:conf/itsc/TuncaliPF16}.

In this paper, we present a specification-driven framework for critical scenario identification (CSI) entirely based on open-source software libraries and demonstrate its benefits with an automated emergency break case study. 
 The proposed framework, based on the falsification testing paradigm~\cite{DBLP:conf/hybrid/NghiemSFIGP10}, uses optimization-based methods for finding critical scenarios.
We first describe the vanilla workflow and show how to tailor it with custom test generation and monitoring strategies. 
 Hence, our aim is to share our experience in combining existing methods into a flexible and efficient SBT framework.
To innovate the methodology for SBT within the framework, 
we investigate the separation between the AD system and the other road users, modeling their interplay with Assume/Guarantee (A/G) contracts. By using A/G contracts, we can improve the search for meaningful scenarios, assign responsibility for critical situations and distinguish between invalid behaviors originating from the AD system and from its environment. In that way, we can detect the violation of environment assumptions in the simulation execution, and discard the test run.
By sharing our experience in SBT, we intend to nurture the innovation of prospective CSI methods that are based on specification-guided strategies.

\section{Specification-Driven Scenario-Based Testing}\label{sec:scenario}
\subsection{Traffic scenario description}

\begin{figure}[tb]
    \centering
    \includegraphics[width=0.70\linewidth]{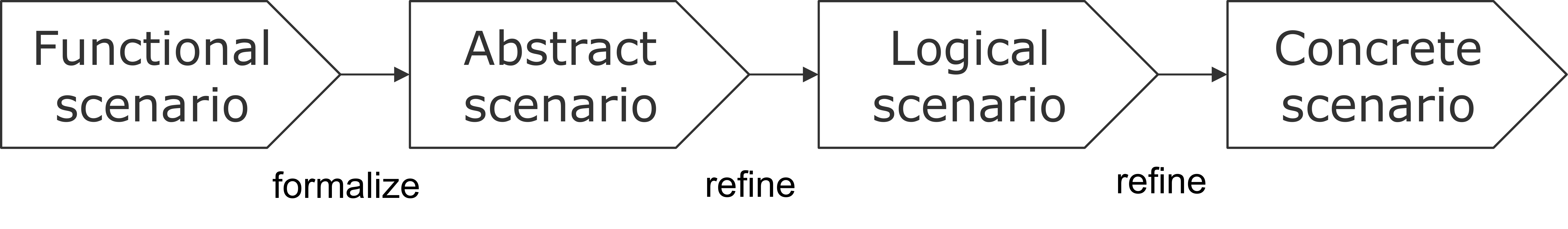}
    \caption{Scenario abstraction types according to \cite{neurohr2021criticality}.}
    \label{fig:abstraction_types}
\end{figure}

In the operational domain in which the ADS will be deployed, it is exposed to a potentially infinite number of traffic scenarios.
As a consequence, it is impractical to conduct testing - even in simulation - directly on these traffic scenarios. A first step towards a successful application of scenario-based testing to assure the correct behavior of an ADS within its ODD is the abstraction of traffic scenarios. While the argumentation for quality assurance is done on a higher level of  abstraction, the creation of evidence is performed on simulating a variety of concrete traffic scenarios derived from the abstract ones. The PEGASUS project ``for the establishment of generally accepted quality criteria, tools and methods as well as scenarios and situations for the release of highly-automated driving functions", introduced three abstraction types: functional, logical, concrete scenarios \cite{menzel2018scenarios}. In this paper, we use an extended classification proposed in \cite{neurohr2021criticality}, see Fig.~\ref{fig:abstraction_types}. Functional scenarios are defined as behavior-based, non-formal descriptions of traffic scenarios in natural language. Abstract scenarios are a formalization of functional scenarios using a declarative way to describe the scenario. Logical scenarios are defined as a parameterized set of traffic scenarios, while concrete scenarios are instances of a logical scenario with fixed parameters. They have a fixed scenery and road user behavior, that is based on the ego-vehicle movement.
Abstract, logical, and concrete scenarios are machine-readable, and various realizations of traffic scenario description formats exist for simulation.
\begin{table}[tb]
\caption{Supported types and properties of scenario description formats }
\label{tab:scenario_formats}
\centering
\begin{tabular}{p{0.23\textwidth}C{0.25\textwidth}C{0.25\textwidth}C{0.18\textwidth}}\toprule
                                                                           & \textbf{OSC1.2} & \textbf{OSC2.0} & \textbf{SCENIC} \\ \midrule
\begin{tabular}{C{0.24\textwidth}}
\textbf{Scenario types}\\ \midrule \end{tabular} & & & \\
Functional                                                         & $\times$                &    $\times$              &    $\times$     \\
Abstract                                                           & $\times$                  &       \checkmark          &    \checkmark    \\
Logical                                                            & \checkmark                &        \checkmark         &    \checkmark    \\
Concrete                                                           & \checkmark                 &  \checkmark               &    \checkmark    \\ \midrule
\begin{tabular}{C{0.24\textwidth}}
\textbf{Properties}\\ \midrule \end{tabular} & & & \\
Syntax format                                                              &    XML\tablefootnote{Extensible Markup Language}             &  DSL\tablefootnote{Domain-specific language}, pythonic               &    DSL, pythonic    \\
Language paradigm                                     &            imperative         &         mostly declarative        & declarative/ imperative     \\
\begin{tabular}[c]{@{}l@{}}Map-agnostic\\ scenario definition\end{tabular} &        $\times$         &       \checkmark          &    \checkmark    \\ \bottomrule
\end{tabular}
\end{table}
In the following, we give a comparison between three non-proprietary, and openly available scenario description formats: OpenSCENARIO\textregistered1.2 (OSC1.2) \cite{asam_openscenario_v1}, OpenSCENARIO\textregistered2.0 (OSC2.0) \cite{asam_openscenario_v2}, and Scenic \cite{fremont2019scenic}, see Tab.~\ref{tab:scenario_formats}. With regard to the overall traffic scenario, their focus is on the initial placement and the dynamic behavior of the actors. The description of the scenery, such as the map, is defined outside these formats.
OSC1.2 is mainly used for describing concrete traffic scenarios that can be directly run by the simulator. The actors' placement and behavior are defined in an imperative fashion using pairs of actions and triggers that evoke these actions. OSC2.0's and Scenic's main intent is to define abstract scenarios, which can be concretized by a dedicated scenario generation engine. OSC2.0's description is mostly declarative by constraining the road users' behavior. The probabilistic programming language Scenic is declarative in the initial actor placement with a rich instruction set for relationships between entities, and uses an imperative description for behaviors. All three languages support parameterization of scenario parameters to describe logical scenarios.
A distinctive feature of OSC2.0 and Scenic compared to OSC1.2 is that the location of the scenario does not need to be specified within the scenario definition. Instead, the scenario generation engine will find a suitable segment on the road map, on which the scenario can be executed with all actors in the simulator.

Based on the scenario format, a database of abstract/logical scenarios needs to be created that covers all the relevant features in the considered ODD of the AD function.
In this paper, we selected Scenic as our scenario format, due to both its flexibility in expressing abstract scenarios and the availability of an open-source testing framework \cite{verifai-cav19} that is provided for Scenic. 

\subsection{Critical scenario identification}
\begin{figure}[tb]
    \centering
    \includegraphics[width=0.98\linewidth]{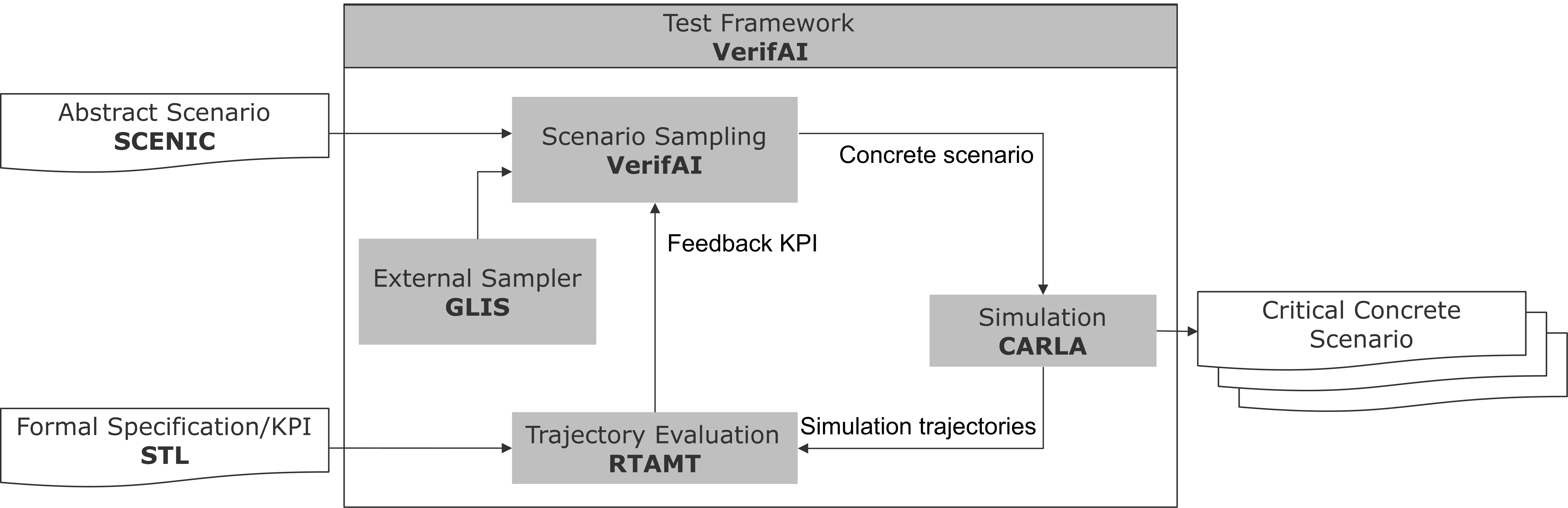}
    \caption{Critical scenario identification framework with tool architecture.}
    \label{fig:csi_framework}
\end{figure}

This section introduces the test framework to find critical concrete scenario instances within a specified abstract scenario efficiently and in a flexible manner. The framework depicted in Fig.~\ref{fig:csi_framework} indicating the overall workflow is based on open-source software components highlighted in bold. It assumes two inputs, the abstract scenario given in the Scenic format, and a formal specification of the AD system defined in signal temporal logic (STL), that we use as a test oracle. 
The technical details on the formal specification are introduced in Sec.~\ref{sec:formal_spec}. 

\paragraph{Workflow} The test execution framework is based on Berkeley's VerifAI \cite{verifai-cav19}. By applying a sampling strategy, VerifAI generates concrete scenarios from the Scenic scenario that are executed in the CARLA simulator \cite{Dosovitskiy17}.
To evaluate the resulting trajectories, we integrated RTAMT - an STL monitoring library \cite{nivckovic2020rtamt} - into the VerifAI-based testing framework. RTAMT provides the automated generation of robustness monitors from STL specifications and therefore facilitates checking simulation traces against the formal specification. The robustness measure is then fed back as a criticality indicator to the scenario sampler that determines new test parameters that constitute the next concrete scenario to be simulated. Depending on the sampling strategy, the scenario search can be of explorative or exploitative nature. Instead of using the sampling strategies provided by VerifAI, we integrated an external sampling strategy, that is based on the global optimizer GLIS \cite{Bem20}.
The details about GLIS are outlined in Sec.~\ref{sec:glis}.

\subsection{Formal Specifications}\label{sec:formal_spec}

Concrete scenarios are typically evaluated against requirements. These requirements can cover various aspects, including safety, legal, comfort and ethical considerations. In order to avoid ambiguities and facilitate their evaluation, there is a need to formulate requirements using a formal specification language. In this paper, we adopt signal temporal logic (STL)~\cite{DBLP:conf/formats/MalerN04} as our specification formalism. There are several motivations to choose STL for requirement formalization: (1) an existing body of work already captures AD system requirements using STL, (2) STL admits quantitative semantics that can be used to guide the search for critical scenarios, and (3) there are runtime verification tools that enable evaluation of STL properties. 
The syntax of STL is given by the grammar
\begin{align*}
\f ::= \true \mid f(\ares) > 0 \mid \lnot \f \mid \f_1 \lor \f_2 \mid \f_1 \until_I \f_2 \mid \f_1 \since_I \f2 \,,
\end{align*}
where $f(R)$ are terms in $\tees$ and $I$ are real intervals with bounds in $\mathbb Q_{\ge 0} \cup \{\infty\}$.
As customary we use $\eventually_I \f \equiv \true \until_I \f$ for \emph{eventually}, $\always_I \f \equiv \lnot \eventually_I \lnot \f$ for \emph{always}, $\once_I \f \equiv \true \since_I \f$ for \emph{once} and $\historically_I \f \equiv \lnot \once_I \lnot \f$ for \emph{historically}.
The timing interval $I$ may be omitted when $I=[0,\infty)$ or $I=(0,\infty)$. STL can be naturally equipped with \emph{quantitative semantics} based on the infinity norm~\cite{DBLP:conf/formats/DonzeM10} that measure how far is the observed behavior from satisfying or violating a requirement.



The evaluation of an AD system cannot be performed in isolation from its environment. For instance, an AD system cannot guarantee safety requirements, such as RSS, in presence of other road users that do not behave in a reasonable manner. 
The relation of the AD system and the environment under which it operates can be formalized in terms of a \emph{contract} $C = (\phi, \psi)$, a pair of properties where $\phi$ represents the assumptions on the environment and $\psi$ guarantees of the system under these assumptions. This classical interpretation of $C$ is given by the temporal logic formula
$$
\always \phi \to \always \psi.
$$
According to the above formula, any violation of the assumption by the environment results in the (vacuous) satisfaction of the contract, even if the system also violates its guarantee. However, this definition neglects that these two violations may not be causally related -- the violation of $\psi$ by the system at time $t$ before the violation of $\phi$ by the environment at time $t' > t$ still results in the satisfaction of the contract. To address this situation, we propose a more refined notion of a contract that takes the intended temporal causality between the environment and the AD system into account. We denote our refined contract by $\hat{C}$ and capture its meaning using the formula:  
$$
\always( (\historically_{[0,T]} \phi) \to \psi)
$$
\noindent where $T$ specifies the maximum duration within which we consider the violation of $\phi$ to be causally related to the violation of $\psi$.

\subsection{Sampling Strategy}\label{sec:glis}
Different sampling strategies may be used to identify the parameters of the next concrete scenario to simulate. These strategies can be broadly divided into na\"ive (passive) and guided search (active) sampling strategies~\cite{zhang2022finding}. The na\"ive search strategies, such as random sampling, involve the independent selection of test parameters. In contrast, the guided search, such as optimization~\cite{feng2020testing_1,feng2020testing_2}, make the selection based on a specific selection criterion and the information of existing samples. Na\"ive search sampling strategies are useful if the simulation is computationally cheap to run since parallelization of the procedure is possible due to the independence among testing samples.
On the other hand, when the test case simulation is computationally expensive to run and/or when the test cases interested (critical test cases in this case) are in a small region of the search domain, the guided search sampling strategies can be more sample efficient.

For the current study, guided-search sampling strategies such as surrogate-based black-box optimization methods are appropriate to efficiently identify relevant critical concrete scenarios for the AD system. It is because a closed-form expression of the KPI in terms of the test parameters is often unavailable. 
Specifically, we use the global optimization algorithm GLIS (Global optimization via Inverse distance weighting and Surrogate radial basis functions)~\cite{Bem20} as the active guided-search sampler to identify the next test parameters of a concrete scenario for testing.
The procedure of GLIS includes an initial sampling stage and an active learning stage. In the initial sampling stage, $N_{\rm initial}$ different test parameters are randomly selected within the search domain, and the corresponding concrete scenarios are simulated. The resulting quantitative evaluation of each test parameter from RTAMT monitors is fed back to GLIS (c.f. Fig.~\ref{fig:csi_framework}). A surrogate radial basis interpolation function (RBF) representing the correlation between the test parameters and the KPI is fitted to the initial samples. In the active learning phase, at each iteration, we identify a new test parameter, simulate the corresponding concrete scenario, and refit the surrogate function by including the newly identified test parameter and its KPI. The new test parameter is obtained by optimizing an acquisition function, which trades off the exploitation of the fitted RBF surrogate and exploration of an inverse distance weighting (IDW) function. IDW is a distance-based exploration function that promotes visiting points far away from the existing samples, which helps prevent the solver from being trapped in the local optima. GLIS terminates when the maximum allowed iteration is reached, or another user-defined criterion is met.

GLIS is chosen for this study, as it easily incorporates constraints and has a low computing cost~\cite{Bem20}. If the computing cost is reasonable, GLIS may be replaced by other surrogate-based active samplers, such as Bayesian optimization.

\section{Automatic Emergency Braking Case Study}\label{sec:simulation}
\begin{figure}[tb]%
    \centering
    \subfloat{ \includegraphics[width=0.4\linewidth, trim={6.2cm 0.01cm 6.2cm 1cm}, clip]{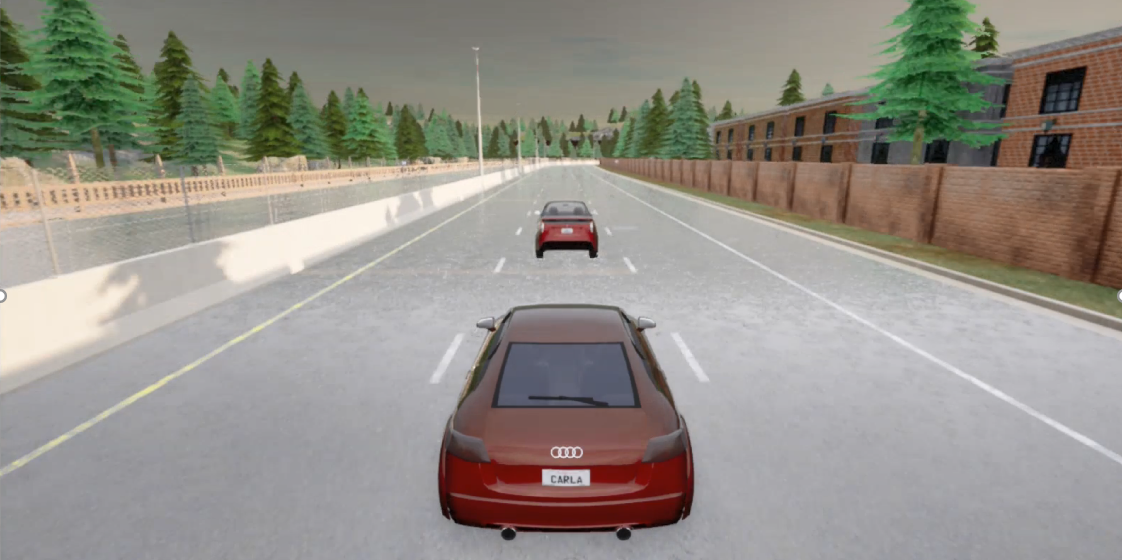} }%
    \qquad
    \subfloat{\includegraphics[width=.52\linewidth, trim={2.5cm 1.5cm 2.5cm 2.5cm},clip]{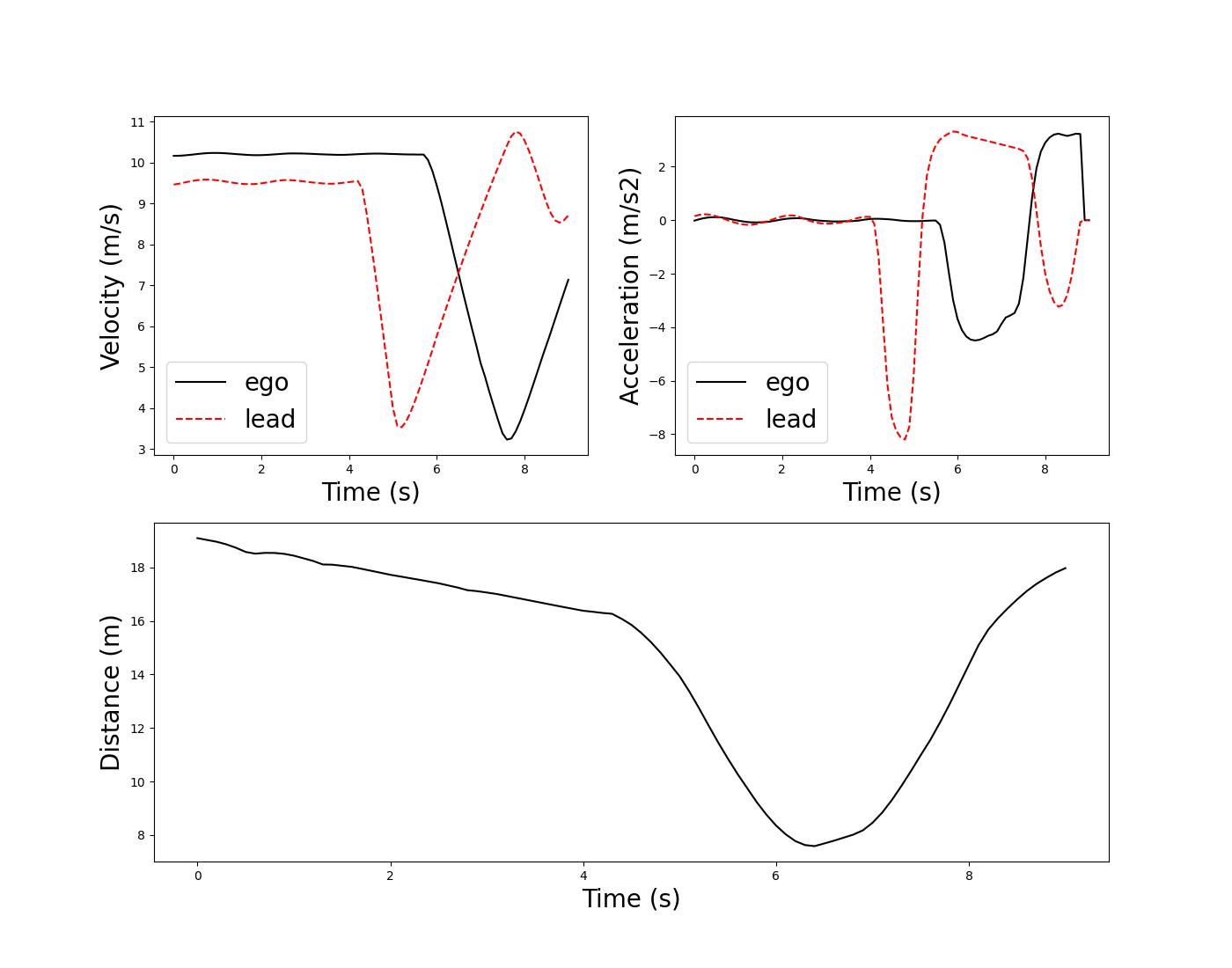} }%
     \setlength{\belowcaptionskip}{-12pt} 
    \caption{(left) Snapshot of CARLA simulator running the AEB function test on a highway. (right) Example of telemetric data collected from all actors.}%
    \label{fig:highway_carla}%
\end{figure}

To illustrate the methodology, we focus on testing a simple Automatic Emergency Braking (AEB) functionality using a highway scenario. 
\paragraph{Scenario description} The functional scenario is an ego vehicle following a leading vehicle on a highway, when suddenly the leading vehicle brakes abruptly. The ego vehicle is equipped with a simplistic distance-based AEB function which is activated when the ego is less than \texttt{safeDist} meters from the leading vehicle. Figure \ref{fig:highway_carla} shows a snapshot of the scenario running in CARLA v9.10.


The abstract scenario, depicted in Listing~\ref{fig:scenic_snippet}, is formulated using Scenic\footnote{The shown scenario is simplified to facilitate presentation.}. The scenario first specifies the sampler and the map used to generate concrete simulations (lines 1 and 2). Then, it defines parameter variables that we partition into: (1) the \emph{constant} variables (lines 3-4) that do not change across concrete scenarios and (2) the \emph{optimization} 
variables (lines 6-7) that are fed to an external (VerifAI) sampler in order to find critical scenarios in a controlled fashion. There are also what we call \emph{implicit} variables that are not explicitly part of the Scenic abstract scenario but still need to have a concrete value in the simulator. For example, the weather conditions, the exact starting position and orientation of each vehicle, the vehicle model, etc. In this case study, there are more than $25$ implicit parameters. The scenario also defines the behavior of the ego (lines 9-12) and of the lead vehicle (lines 14-18). Both the ego and the lead vehicle follow the lane with some target speed as their default behavior. However, the lead vehicle abruptly breaks at regular intervals, while the ego breaks when it approaches any object at some minimum distance. The two vehicles are spawned at some uniformly chosen part of the map (line 20) that is sufficiently far away from an intersection (line 25). The lead car is initialized at some pre-defined distance in front of the ego vehicle (lines 22-23).



\begin{lstlisting}[
    basicstyle=\ttfamily\footnotesize,
    language=SCENIC,
    frame=tb,
    caption= AEB highway scenario in Scenic, 
    captionpos=b,
    belowcaptionskip=2em,
    label=fig:scenic_snippet
]
param verifaiSamplerType = `glis' # specify sampler 
param carla_map = `Town04' # specify map to use
initDist = 30 # constant 
leadSpeed = 10

safeDist = VerifaiRange(25,45) # optimization variable 
egoSpeed = VerifaiRange(9,11) # optimization variable 

behavior AEB_Behavior: # define ego behavior
  try:
    FollowLaneBehavior(egoSpeed)
  interrupt when withinDistanceToObjsInLane(self, safeDist):   take SetBrakeAction(1), SetThrottleAction(0)

behavior Brake_Behavior: # define behavior of lead car
  try:
    FollowLaneBehavior(leadSpeed)
  interrupt when simulation().currentTime > delay:
    take SetBrakeAction(1), SetThrottleAction(0)
    
spawnPt = Uniform(*HighwayRoads) # Highway part of map

ego = Car at spawnPt, with AEB_Behavior} # spawn ego  
leadCar = Car at spawnPt + initDist, with Brake_Behavior # lead

require (distance from leadCar to intersection > 50)
# extra requirements for rejection sampling
\end{lstlisting}

\paragraph{Formalized Requirements}

We illustrate the formalization of the requirements with the contract $C=(\varphi, \psi)$, which captures the assumption $\varphi$ about the maximum allowed deceleration of the lead vehicle and the guarantee $\psi$ as the Responsibility-Sensitive Safety (RSS) property of the ego vehicle.
The assumption $\varphi$ originates from the IEEE Standard 2846-2022~\cite{ieee2846}, that describes the minimal set of assumptions on the road users for safety-related models of AD. From the assumptions described in the standard, we focus on the maximum deceleration specification
$$
\varphi = \beta \leq \beta_{\emph{max}}.
$$

The Responsibility-Sensitive Safety (RSS) rule specifies, under minimal assumptions, what longitudinal and lateral distances the ego vehicle must keep from other road users to ensure no collisions \cite{ShalevSchwartz2017}. The RSS rules were formalized into temporal logic by \cite{Arechiga2019} \cite{Hekmatnejad2019}. We adopt the STL specification from~\cite{Arechiga2019} for an ego vehicle (\textit{back}) to keep a safe longitudinal distance to another vehicle (\textit{front}):
\begin{align*}
    &\always \left( v_{\text{front}} \geq 0 \wedge v_{\text{back}} \geq 0 \right) \\
    &\always \left( a_{\text{front}} \in [a_{\text{max-Br}}, a_{\text{max-Acc}}] \wedge a_{\text{back}} \in [a_{\text{max-Br}}, a_{\text{max-Acc}}] \right) \\
    &\always \left( d(\text{front, back}) < d_{\text{safe}} \rightarrow a_{\text{back}} \in [a_{\text{max-Br}}, a_{\text{min-Br}}] \right) 
\end{align*}
where $a,v$ are correspondingly acceleration and velocity. Similarly $a_{\text{max-Acc}}$, $a_{\text{max-Br}}$, $a_{\text{min-Br}}$ are assumed maximum acceleration, maximum braking, and minimum braking acceleration. Finally, $d_{\text{safe}}$ is determined dynamically depending on the velocities of both vehicles, and the reaction time $\tau$ of the ego vehicle:
\begin{equation*}
    d_{\text{safe}} = \left( v_{\text{back}} \tau + \frac{a_{\text{max-Acc}} \tau^2}{2} + \frac{(v_{\text{back}} + a_{\text{max-Acc}} \tau )^2}{2 a_{\text{min-Br}}} - \frac{v_{\text{front}}^2}{2 a_{\text{max-Br}}} \right).
\end{equation*}
The safety distance is calculated in order to ensure that a collision is avoided as long as the ego vehicle is sufficiently far away from the leading vehicle. If it is momentarily closer than $d_{\text{safe}}$ then a collision will still be avoided if the ego is reacting appropriately (by braking with at least $a_{\text{min-Br}}$).

\subsection{Simulation Results}
\begin{figure}[htb]%
    \centering
    \subfloat{{\includegraphics[width=.46\linewidth,clip]{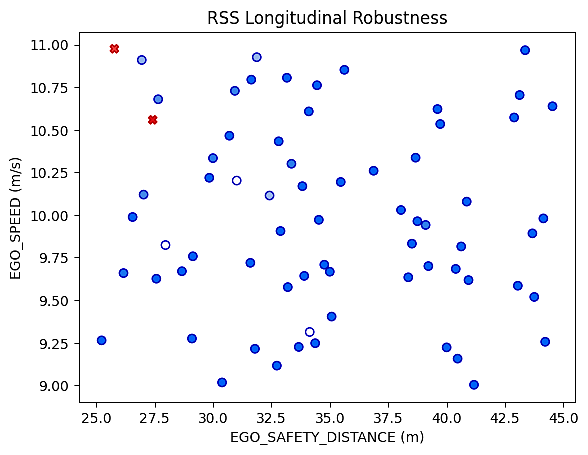} }}%
    \qquad
    \subfloat{{\includegraphics[width=.46\linewidth,clip]{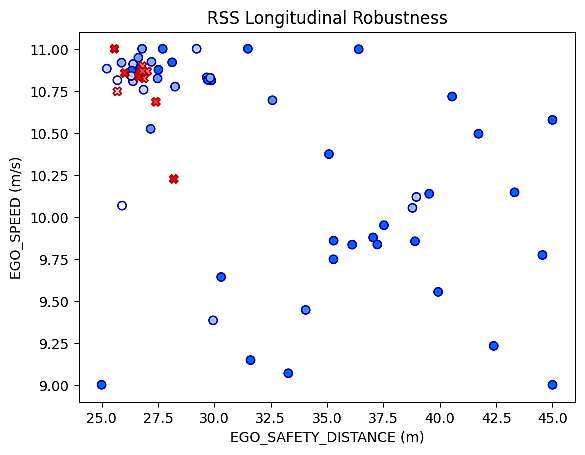} }}%
     \setlength{\belowcaptionskip}{-12pt} 
    \caption{ Comparison between Halton sampling (left) and GLIS  sampling \cite{Bem20}(right) for $70$ concrete scenarios. The GLIS parameters are: $\alpha=1$, $\delta=0.5$, $\varepsilon_{\text{SVD}}=0.01$, and an inverse-quadratic basis function with $\epsilon=0.2$ was used.}%
    \label{fig:simulation_results}%
\end{figure}

In this section we present our evaluation outcomes.
Fig. \ref{fig:simulation_results} shows the results from simulating the abstract scenario $70$ times using the described tool chain. Each point in the scatter plots represents a simulated concrete scenario, where the RSS longitudinal distance was monitored. If the ego vehicle managed to react adequately by braking in time, then this is represented as a blue circle, otherwise (if the specification was violated) it is represented by a red cross (the intensity of the color represents the robustness degree). 

Furthermore, we compare two different sampling strategies to find critical scenarios. In this case, we compare a passive sampling strategy (i.e.\ agnostic to feedback) which is based on Halton sequences~\cite{DBLP:journals/cacm/Halton64}, to an active strategy based on the GLIS optimization sampling. As expected, sampling scenarios with GLIS leads to the discovery of more critical scenarios ($11$ compared to $2$ with Halton), and suggests variable regions which should be further investigated. In our example, the optimizer clearly was trying to exploit around the region of higher \texttt{egoSpeed}, and lower \texttt{safeDist} (as expected). In practice, both strategies are used to obtain a clear picture of the performance of the ADAS functionality.

\begin{figure}[htb]
    \centering
    \includegraphics[width=\linewidth]{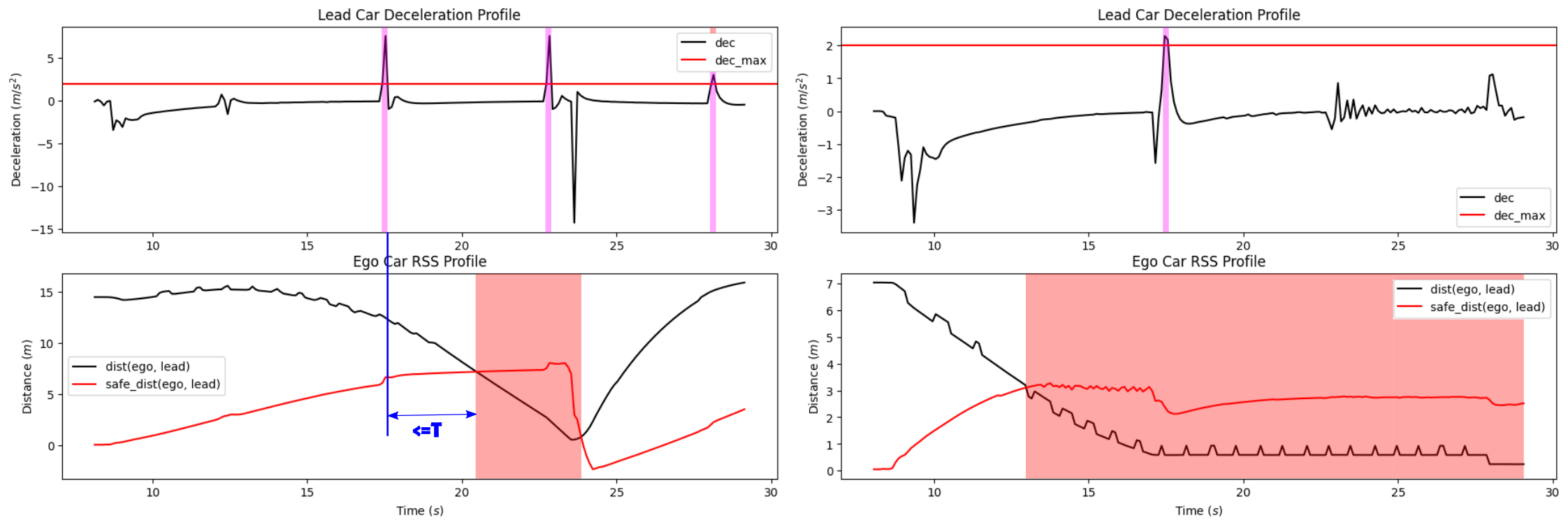}
    \caption{Evaluation with A/G contracts.}
    \label{fig:ag}
\end{figure}

In Fig.~\ref{fig:ag}, we illustrate the discrepancy between the classical and the refined interpretation of A/G contracts. The figure depicts two simulations showing the deceleration $\beta$ of the lead vehicle and the maximum allowed deceleration threshold $\beta_{\emph{max}}=2 m/s^2$ (top) and the distance between the ego and the lead vehicle, as well as the safe distance between them (bottom). We see that in the two simulations both the assumption $\varphi$ and the guarantee $\psi$ are violated (purple and red stipes, respectively). In the first simulation (left), there is a clear causality between the abrupt breaking of the lead vehicle and the longitudinal RSS violation -- it follows that the contract is satisfied under both the classical and the refined interpretation. In the second simulation (right), the violation of the longitudinal RSS requirement happens before the lead vehicle breaks. Intuitively, we expect the contract to be violated since the behavior of the lead vehicle did not cause this critical scenario. However, under the classical contract interpretation, the contract is satisfied because the lead vehicle does violate the assumption at a later stage. On the other hand, the refined contract rightly indicates the contract falsification.

\subsection{Lessons Learned}
\label{sec:lessons}

In this section, we share our experience about the scenario-based testing framework and collected during the case study evaluation.

\paragraph{Passive vs. active sampling} Both passive and active sampling have their merits in testing AD systems. Passive sampling methods such as Halton provide a coverage of the parameter space, facilitate detecting interesting patterns, if any, and help identifying parameter regions that are interesting to further explore. In contrast, active sampling methods such as GLIS can accelerate the detection of critical scenarios.

\paragraph{Level of scenario abstraction} Balance between keeping a scenario abstract, and letting the tools sample different variables, and having consistent concrete scenarios. If too many variables are left unspecified, drawing meaningful conclusions from the experiments is difficult, but if too many parameters are specified, there is a risk of missing out on potential critical scenarios that are relevant (and it also needs more development time). 

\paragraph{Optimization with implicit variables} It is interesting to note that from the point of view of the optimizer, the robustness function the of concrete scenario is non-deterministic. That is, there are many different concrete scenarios that result from having the same \texttt{egoSpeed} and \texttt{safeDist} which result in different robustness values. This is mostly due to different implicit parameters impacting the robustness, which the optimizer does not directly see (e.g.\ road geometry).


\bibliographystyle{splncs04}
\bibliography{bibliography}

\end{document}

%% file: macros.tex
\let\phi\varphi

\newcommand\f\varphi
\newcommand\g\psi

\newcommand\tees{\Theta}

\font\zmsy=zmsy10
\newcommand\diaM{\hbox{\zmsy\char'121}}
\newcommand\boxM{\hbox{\zmsy\char'140}}

\DeclareMathOperator\once{\diaM}
\DeclareMathOperator\historically{\boxM}
\DeclareMathOperator\eventually{\lozenge}
\DeclareMathOperator\always{\Box}

\DeclareMathOperator\until{\mathcal{U}}
\DeclareMathOperator\since{\mathcal{S}}
\newcommand\true{\top}

\newcommand\ares{R}



%% file: main.bbl
\begin{thebibliography}{10}
\providecommand{\url}[1]{\texttt{#1}}
\providecommand{\urlprefix}{URL }
\providecommand{\doi}[1]{https://doi.org/#1}

\bibitem{ieee2846}
{IEEE Standard for Assumptions in Safety-Related Models for Automated Driving
  Systems}. IEEE Std 2846-2022 pp. 1--59 (2022).
  \doi{10.1109/IEEESTD.2022.9761121}

\bibitem{Arechiga2019}
Aréchiga, N.: Specifying safety of autonomous vehicles in signal temporal
  logic. In: 2019 IEEE Intelligent Vehicles Symposium (IV). pp. 58--63 (2019).
  \doi{10.1109/IVS.2019.8813875}

\bibitem{asam_openscenario_v1}
{Association for Standardization of Automation and Measuring Systems}: {ASAM
  OpenSCENARIO V1.2.0}. Standard, Munich, Germany (2022),
  \url{https://www.asam.net/standards/detail/openscenario/}

\bibitem{asam_openscenario_v2}
{Association for Standardization of Automation and Measuring Systems}: {ASAM
  OpenSCENARIO V2.0.0}. Standard, Munich, Germany (2022),
  \url{https://www.asam.net/standards/detail/openscenario/v200/}

\bibitem{Bem20}
Bemporad, A.: Global optimization via inverse distance weighting and radial
  basis functions. Computational Optimization and Applications  \textbf{77},
  571--595 (2020), code available at
  \url{http://cse.lab.imtlucca.it/~bemporad/glis}

\bibitem{DBLP:conf/formats/DonzeM10}
Donz{\'{e}}, A., Maler, O.: Robust satisfaction of temporal logic over
  real-valued signals. In: Formal Modeling and Analysis of Timed Systems - 8th
  International Conference, {FORMATS} 2010, Klosterneuburg, Austria, September
  8-10, 2010. Proceedings. pp. 92--106 (2010)

\bibitem{Dosovitskiy17}
Dosovitskiy, A., Ros, G., Codevilla, F., Lopez, A., Koltun, V.: {CARLA}: {An}
  open urban driving simulator. In: Proceedings of the 1st Annual Conference on
  Robot Learning. pp. 1--16 (2017)

\bibitem{DBLP:journals/jar/DreossiDS19}
Dreossi, T., Donz{\'{e}}, A., Seshia, S.A.: Compositional falsification of
  cyber-physical systems with machine learning components. J. Autom. Reason.
  \textbf{63}(4),  1031--1053 (2019). \doi{10.1007/s10817-018-09509-5},
  \url{https://doi.org/10.1007/s10817-018-09509-5}

\bibitem{verifai-cav19}
Dreossi, T., Fremont, D.J., Ghosh, S., Kim, E., Ravanbakhsh, H.,
  Vazquez{-}Chanlatte, M., Seshia, S.A.: {VerifAI:} {A} toolkit for the formal
  design and analysis of artificial intelligence-based systems. In: 31st
  International Conference on Computer Aided Verification (CAV) (Jul 2019)

\bibitem{feng2020testing_2}
Feng, S., Feng, Y., Sun, H., Bao, S., Zhang, Y., Liu, H.X.: Testing scenario
  library generation for connected and automated vehicles, part ii: Case
  studies. IEEE Transactions on Intelligent Transportation Systems
  \textbf{22}(9),  5635--5647 (2020)

\bibitem{feng2020testing_1}
Feng, S., Feng, Y., Yu, C., Zhang, Y., Liu, H.X.: Testing scenario library
  generation for connected and automated vehicles, part i: Methodology. IEEE
  Transactions on Intelligent Transportation Systems  \textbf{22}(3),
  1573--1582 (2020)

\bibitem{fremont2019scenic}
Fremont, D.J., Dreossi, T., Ghosh, S., Yue, X., Sangiovanni-Vincentelli, A.L.,
  Seshia, S.A.: Scenic: a language for scenario specification and scene
  generation. In: Proceedings of the 40th ACM SIGPLAN Conference on Programming
  Language Design and Implementation. pp. 63--78 (2019)

\bibitem{DBLP:journals/cacm/Halton64}
Halton, J.H., Smith, G.B.: Algorithm 247: Radical-inverse quasi-random point
  sequence. Commun. {ACM}  \textbf{7}(12),  701--702 (1964).
  \doi{10.1145/355588.365104}, \url{https://doi.org/10.1145/355588.365104}

\bibitem{Hekmatnejad2019}
Hekmatnejad, M., Yaghoubi, S., Dokhanchi, A., Amor, H.B., Shrivastava, A.,
  Karam, L., Fainekos, G.: Encoding and monitoring responsibility sensitive
  safety rules for automated vehicles in signal temporal logic. In: Proceedings
  of the 17th ACM-IEEE International Conference on Formal Methods and Models
  for System Design. MEMOCODE '19 (2019). \doi{10.1145/3359986.3361203},
  \url{https://doi.org/10.1145/3359986.3361203}

\bibitem{DBLP:conf/formats/MalerN04}
Maler, O., Nickovic, D.: Monitoring temporal properties of continuous signals.
  In: Formal Techniques, Modelling and Analysis of Timed and Fault-Tolerant
  Systems, Joint International Conferences on Formal Modelling and Analysis of
  Timed Systems, {FORMATS} 2004 and Formal Techniques in Real-Time and
  Fault-Tolerant Systems, {FTRTFT} 2004, Grenoble, France, September 22-24,
  2004, Proceedings. pp. 152--166 (2004)

\bibitem{menzel2018scenarios}
Menzel, T., Bagschik, G., Maurer, M.: Scenarios for development, test and
  validation of automated vehicles. In: 2018 IEEE Intelligent Vehicles
  Symposium (IV). pp. 1821--1827. IEEE (2018)

\bibitem{neurohr2021criticality}
Neurohr, C., Westhofen, L., Butz, M., Bollmann, M.H., Eberle, U., Galbas, R.:
  Criticality analysis for the verification and validation of automated
  vehicles. IEEE Access  \textbf{9},  18016--18041 (2021)

\bibitem{DBLP:conf/hybrid/NghiemSFIGP10}
Nghiem, T., Sankaranarayanan, S., Fainekos, G., Ivancic, F., Gupta, A., Pappas,
  G.J.: Monte-carlo techniques for falsification of temporal properties of
  non-linear hybrid systems. In: Proceedings of the 13th {ACM} International
  Conference on Hybrid Systems: Computation and Control, {HSCC} 2010,
  Stockholm, Sweden, April 12-15, 2010. pp. 211--220 (2010)

\bibitem{nivckovic2020rtamt}
Ni{\v{c}}kovi{\'c}, D., Yamaguchi, T.: Rtamt: Online robustness monitors from
  stl. In: International Symposium on Automated Technology for Verification and
  Analysis. pp. 564--571. Springer (2020)

\bibitem{orad_j3016}
{On-Road Automated Driving (ORAD) committee}: {J3016 Taxonomy and Definitions
  for Terms Related to Driving Automation Systems for On-Road Motor Vehicles}.
  Tech. rep. (2021),
  \url{{https://www.sae.org/standards/content/j3016\_202104/}}

\bibitem{qin2019automatic}
Qin, X., Ar{\'e}chiga, N., Best, A., Deshmukh, J.: Automatic testing with
  reusable adversarial agents. arXiv preprint arXiv:1910.13645  (2019)

\bibitem{ShalevSchwartz2017}
{Shalev-Shwartz}, S., {Shammah}, S., {Shashua}, A.: {On a Formal Model of Safe
  and Scalable Self-driving Cars}. arXiv e-prints arXiv:1708.06374 (Aug 2017)

\bibitem{tuncali2019requirements}
Tuncali, C.E., Fainekos, G., Prokhorov, D., Ito, H., Kapinski, J.:
  Requirements-driven test generation for autonomous vehicles with machine
  learning components. IEEE Transactions on Intelligent Vehicles
  \textbf{5}(2),  265--280 (2019)

\bibitem{DBLP:conf/itsc/TuncaliPF16}
Tuncali, C.E., Pavlic, T.P., Fainekos, G.: Utilizing s-taliro as an automatic
  test generation framework for autonomous vehicles. In: 19th {IEEE}
  International Conference on Intelligent Transportation Systems, {ITSC} 2016,
  Rio de Janeiro, Brazil, November 1-4, 2016. pp. 1470--1475 (2016)

\bibitem{zhang2022finding}
Zhang, X., Tao, J., Tan, K., T{\"o}rngren, M., Sanchez, J.M.G., Ramli, M.R.,
  Tao, X., Gyllenhammar, M., Wotawa, F., Mohan, N., et~al.: Finding critical
  scenarios for automated driving systems: A systematic mapping study. IEEE
  Transactions on Software Engineering  (2022)

\end{thebibliography}
